\documentclass[11pt]{elsarticle}
\usepackage{pdfsync}
\usepackage{natbib}
\usepackage{geometry}

\usepackage{amsthm}
\usepackage{amsfonts}
\usepackage{MnSymbol}

\usepackage[curve,tips]{xypic}
\SelectTips{eu}{11}
\UseTips

\newcommand{\colim}{\varinjlim}

\renewcommand{\epsilon}{\varepsilon}

\newcommand{\iso}{\cong}

\newcommand{\normal}{\lhd}
\newcommand{\supp}{\operatorname{supp}}

\newcommand{\clh}{{\scriptstyle\bf H}}

\newcommand{\hf}{\clh\mathfrak F}

\newcommand{\Z}{\mathbb Z}
\newcommand{\N}{\mathbb N}

\newcommand{\fpinfty}{{\FP}_{\infty}}

\newcommand{\fpn}{{\FP}_{n}}

\newcommand{\fn}{{\FFF}_{n}}

\newcommand{\ft}{{\FFF}_{2}}

\newcommand{\FFF}{\operatorname{F}}
\newcommand{\FP}{\operatorname{FP}}

\newcommand{\mono}{\rightarrowtail}
\newcommand{\epi}{\twoheadrightarrow}

\renewcommand{\implies}{\Rightarrow}

\renewcommand{\wr}{\nwsebipropto}

\newcommand{\tor}{\operatorname{Tor}}
\newcommand{\ext}{\operatorname{Ext}}

\title
{Graph-wreath products and finiteness conditions\tnoteref{accepted}}
\author{P. H. Kropholler\fnref{phk}}
\tnotetext[accepted]{Accepted for publication in J. Pure Appl. Algebra on 12th June 2015}
\fntext[phk]{Partially supported by EPSRC Grant EP/N007328/1.} 
\author{A. Martino}
\address{Mathematical Sciences, University of Southampton, Southampton SO17 1BJ, United Kingdom}
\ead{p.h.kropholler@soton.ac.uk}
\ead{a.martino@soton.ac.uk}

\begin{document}

\begin{abstract}
A notion of \emph{graph-wreath product of groups} is introduced. We obtain sufficient conditions for these products to satisfy the topologically inspired finiteness condition type $\fn$. Under various additional assumptions we show that these conditions are necessary. Our results generalize results of Bartholdi, Cornulier and Kochloukova about wreath products. Graph-wreath products of groups include classical permutational wreath products and semidirect products of right-angled Artin groups by certain groups of automorphisms amongst others.
\end{abstract}


\maketitle

\section{Introduction}

In this paper we introduce a concept of \emph{graph-wreath product of groups} which encompasses the notions of restricted wreath product and of right-angled Artin group as well as other graph products of groups. Our Theorem A  gives sufficient conditions for these groups to satisfy the finiteness conditions type $\FFF_n$, $n\ge1$. Note that the property $\FFF_2$ is equivalent to finite presentability, and in the particular case $n=2$, our result has already been proved for permutational wreath products by Cornulier.
Both the results and the methods of proof in this paper have been influenced by the techniques of \cite{Cornulier}. In Theorems B, C, D we examine necessary conditions for the same finiteness properties and here we have built on homological methods of Bartholdi--Cornulier--Kochloukova \cite{BCK}. The finiteness conditions $\FFF_n$ are discussed in many articles, but we refer to the fundamental paper \cite{BB} of Bestvina and Brady for background material because that paper also discusses a number of other ideas relevant to our paper. In \cite{BCK}, results similar to ours but confined to the wreath product case are developed for the properties $\FP_n$ rather than $\FFF_n$. Again, \cite{BB} provides an excellent reference for the ideas and remains the fundamental source of examples distinguishing type $\FP_n$ and type $\FFF_n$ when $n\ge2$.

\subsection{Immediate Applications} We conclude the paper with some applications to produce new examples. Here is a summary of some of these.
\begin{itemize}
\item
For arbitrarily large $n$ there exists an elementary amenable group of type $\FFF_n$ which has a subgroup isomorphic to the iterated wreath product
$$\underbrace{(\dots((\Z\wreath\Z)\wreath\Z))\wreath\Z\dots)\wreath\Z}_n.$$ Note that by contrast, every elementary amenable group of type
$\FFF_\infty$ is virtually soluble of finite rank and has no sections isomorphic to $A\wreath\Z$ for any non-trivial group $A$. For further information, see \cite{Kropholler84} where it is shown that finitely generated soluble groups with no large wreath product sections have finite rank, and \cite{KMN} for an account of elementary amenable groups of type $\fpinfty$.
\item There exist groups of type $\FFF_\infty$ which contain infinite normal locally finite subgroups.
There exist groups of type $\FFF_\infty$ which contain normal free abelian groups of infinite rank. These examples are again in sharp contrast with results on $\hf$-groups of type $\fpinfty$ which state that there is a bound on the orders of finite subgroups and that torsion-free subgroups have finite cohomological dimension. See \cite{Kropholler93} for an explanation of the class $\hf$ and the results that follows for groups in this class.
\end{itemize}

\subsection{Historical remarks}
 In 1961, Baumslag\footnote{We record with sadness that Gilbert Baumslag  died on October 20, 2014, while we were completing the writing of this paper. His work was influential in the field generally and specifically with regard to our project.} \cite{Baumslag1961} proved that the standard restricted wreath product $A\wreath H$ of two groups $A$ and $H$ is finitely presented if and only if both $A$ and $H$ are finitely presented and either $A$ is trivial or $H$ is finite. One implication here is relatively elementary, namely that $A\wreath H$ is finitely presented if $A$ is finitely presented and $H$ is finite. The converse is more subtle and that was the main focus of Baumslag's attention in \cite{Baumslag1961}.
Following this initial result, there seems to have been little further progress until, in 1997, the celebrated paper \cite{BB} of Bestvina and Brady appeared with the first known examples of group of type $\FP_2$ which are not finitely presented and at this time it therefore became of interest to review many results about finitely presented groups to find out if they were also true of arbitrary groups of type $\FP_2$.
In 1998, Baumslag--Bridson--Gruenberg \cite{bbg} extended Baumslag's earlier result by proving that $A\wreath H$ is of type $\FP_2$ if and only if both $A$ and $H$ are of type $\FP_2$ and either $A=1$ or $H$ is finite.

The subsequent developments appear in the papers already alluded to above, \cite{BCK,Cornulier}. All the results of the four papers \cite{Baumslag1961,bbg,BCK,Cornulier} concern wreath products. But Cornulier's proof uses an intermediate system of groups between the free group that plays a role in a presentation and the direct sum of groups that plays the role of base of a wreath product. We have developed results which apply to all of these intermediate groups as well as the wreath product.

The graph-wreath product unites the concepts of permutational wreath product of two groups and graph product of a family of groups over a graph.  Constructions of this kind have been considered by other authors: see \cite{Leary2005} for an example. Our methods build on ideas of Davis, \cite{Davis2012}.

\subsection{Finiteness conditions}

We are concerned with the finiteness condition \emph{type $\FFF_n$} of a group $G$, meaning that there is an Eilenberg--Mac Lane space with finite $n$-skeleton. This is a property enjoyed by all finite groups but is topologically inspired and so may be called a \emph{homotopical finiteness condition}. We shall also need the finiteness conditions \emph{type $\FP_n$} for groups and modules. The following are some elementary facts about these definitions. We refer the reader to one of the standards texts \cite{Bieriqmw}, \cite{Brown} for this and other background material.
The following lemma may be a helpful reminder of how the type $\FP_n$ translates between groups and modules.

\subsection{Lemma}\label{elem}
\emph{A group $G$ is of type $\FP_n$ if and only if the trivial module $\Z$ is of type $\FP_n$ as a $\Z G$-module. Every group of type $\FFF_n$ is of type $\FP_n$. If $H$ is a subgroup of $G$ then $H$ is of type $\FP_n$ as a group if and only if $\Z\otimes_{\Z H}\Z G$ is of type $\FP_n$ as a $\Z G$-module.}

As is evident from the notation for induced module in this Lemma, we choose to work with right modules unless it is stated \emph{left modules}.

\subsection{Wreath products}


The restricted wreath product $A\wreath H$ of two groups $A$ and $H$ has \emph{base} $B$ the set of functions from $H$ to $A$ with finite support and \emph{head} $H$. It is the semidirect product $B\rtimes H$. More generally, if $H$ has a permutation representation through an action on a set $\Omega$ then the restricted permutational wreath product $A\wreath_\Omega H$ is constructed in the same way with base the set of functions from $\Omega$ to $A$ with finite support. We do not have anything to say about the unrestricted wreath product.

\subsection{Simple graphs, cliques, and flag complexes}

By a \emph{simple graph} we mean a $1$-dimensional simplicial complex.
A \emph{clique} in such a graph consists of a finite set of vertices each pair of which are joined by an edge. For $p\ge1$
a \emph{$p$-clique} is a clique with exactly $p$ distinct vertices. The \emph{flag complex generated} by a simple graph $\Gamma$ is the simplicial complex whose $1$-skeleton coincides with $\Gamma$ and in which every $(p+1)$-clique supports a $p$-simplex whenever $p\ge2$. If $X$ is a simplicial complex then $X$ is called a \emph{flag complex} if and only if it is the flag complex generated by its $1$-skeleton. In general, by a \emph{non-face} of $X$ we mean a set of $p+1$ distinct vertices with $p\ge2$ which do not support a $p$-simplex. A \emph{non-edge} is a pair of distinct vertices which are not joined by an edge. Another way of defining a flag complex is to say that it is a simplicial complex in which each non-face has a non-edge.

\subsection{Group actions on graphs}
Let $H$ be a group. By a \emph{simple $H$-graph} we mean a simple graph on which $H$ acts by graph automorphisms: we allow elements of $H$ to invert edges. An action of $H$ on a simple graph $\Gamma$ induces an action on the flag complex generated by $\Gamma$: note that such an action need not be admissible in the sense of Brown \cite{Brown1987}, but becomes admissible on barycentric subdivision.

The following lemma is fundamental. It translates between the language of graph theory and the language of simplicial complexes.

\subsection{Lemma}\label{cliquestoflags}
\emph{
Let $H$ be a group and let $\Gamma$ be a non-empty simple $H$-graph. Let $L$ be the flag complex spanned by $\Gamma$.
Let $m$ and $n$ be non-negative integers.
Let $\Delta_m$ denote the set of $m$-simplices of $L$ and let $L^m$ denote the $m$-skeleton of $L$, for $m\ge0$. Then the action of $H$ extends naturally to $L$ and the following are equivalent:
\begin{enumerate}
\item
$\Z\Delta_m$ is of type ${\FP}_{n}$ as a $\Z H$-module.
\item
$H$ has finitely many orbits of $(m+1)$-cliques and the stabilizer of each $(m+1)$-clique has type $\FP_{n}$.
\end{enumerate}
}
\begin{proof}
Assume that (i) holds.
Since $n\ge0$ we infer from (i) that $\Z\Delta_m$ is finitely generated as a $\Z H$-module. The decomposition of $\Delta_m$ into $H$-orbits gives rise to a direct sum decomposition of $\Z\Delta_m$ as $\Z H$-module. Hence $\Delta_m$ falls into finitely many orbits. An $m$-simplex in the flag complex has $m+1$ vertices which form a clique. So having finitely many orbits of $m$-simplices is the same condition as having finitely many orbits of $(m+1)$-cliques. Choose a set $\Delta_\dagger$ of orbit representatives in $\Delta_m$. Then $\Delta_\dagger$ is finite and $\Z\Delta_m$ is isomorphic to $\bigoplus_{\delta\in\Delta_\dagger}\Z\otimes_{\Z H_\delta}\Z H$. Now $\Z\Delta_m$ is of type $\FP_n$ over $\Z H$ if and only if each $\Z\otimes_{\Z H_\delta}\Z H$ is of type $\FP_n$ over $\Z H$. This in turn is equivalent to each stabilizer $H_\delta$ being of type $\FP_n$ as a group by Lemma \ref{elem}. The converse is equally straightforward and we leave the details.
\end{proof}

\subsection{Graph products}


If $\Gamma$ is a simple graph and $\mathbf A=(A_v)_{v\in V}$ is a family of groups indexed by the vertices of $\Gamma$ then the graph product $\mathbf A^\Gamma$ is defined to be the quotient of the free product formed by imposing commutator relations that force elements over distinct vertices to commute when those vertices are joined by an edge. For example:
\begin{itemize}
\item
If $\Gamma$ is the graph consisting of two vertices and one edge then $(A_1,A_2)^\Gamma$ is the direct product $A_1\times A_2$.
\item
If $\Gamma$ is the graph consisting of two vertices and no edges then $(A_1,A_2)^\Gamma$ is the free product $A_1*A_2$.
\item
If $\Gamma$ is a complete graph then $\mathbf A^\Gamma$ is the direct sum $\bigoplus_{v\in V}A_v$.
\item
If $\Gamma$ is a discrete graph (with no edges) then $\mathbf A^\Gamma$ is the free product $\star_{v\in V}A_v$.
\end{itemize}
If all the vertex groups $A_v$ are equal to the same group $A$ then we write $A^\Gamma$ instead of $\mathbf A^\Gamma$. These further examples belong to this case:
\begin{itemize}
\item
For any graph, $\Z^\Gamma$ is the right-angled Artin group determined by $\Gamma$ and $(\Z/2\Z)^\Gamma$ the right-angled Coxeter group.
\end{itemize}

\section{The graph-wreath product and statement of results}

We now come to the main definition of the paper. Constructions of the following kind have certainly been used frequently when $H$ is finite, but less commonly when $H$ is infinite. Leary has used the construction in case $H$ is an infinite cyclic group, \cite{Leary2005}.

\subsection{Definition}
Suppose that $H$ is a group acting on a graph $\Gamma$. Given a group $A$ we can form the graph product of the family in which the same group $A$ is placed over each vertex of $\Gamma$. Then the action of $H$ on $\Gamma$ induces an action of $H$ on the graph product $A^\Gamma$. We define the \emph{graph-wreath product} to be the semidirect product $A^\Gamma\rtimes H$. We introduce the notation $A\wr_\Gamma H$ for this construction. When $H$ is trivial this is just the graph product. When $\Gamma$ is the complete graph on its vertex set then this is the wreath product. If
$\buildrel\leftrightarrow\over\Omega$ denotes the complete graph on the $H$-set $\Omega$ then
$A\wr_{\buildrel\leftrightarrow\over\Omega}H=A\wreath_\Omega H$.
If $\Gamma$ has two vertices and one edge and $C_2$ is the group of order two that inverts the edge then $A\wr_\Gamma C_2$ is the wreath product $A\wreath C_2$. If on the other hand, $\Gamma$ has two vertices and no edge then $A\wr_\Gamma C_2$ is isomorphic to the free product $A*C_2$. The notation $\wr$ is intended to distinguish our definition from the wreath product while recognising a connection.

Our results fall naturally into two categories: sufficient conditions which are easier to establish and may have wider immediate applicability; and the necessity of those sufficient conditions under various additional hypotheses. In the end it remains open to find a clean set of necessary and sufficient conditions.

\subsection{Main results}\label{mainstatement}

We fix the following notation: $H$ is a group, $\Gamma$ is a non-empty simple $H$-graph and $A$ is another group. We write $V$ for the set of vertices of $\Gamma$ and we write $L$ for the flag complex spanned by $\Gamma$. Note that $V$ is also the set of vertices of $L$. Let $G=A\wr_\Gamma H$ be a graph-wreath product of $A$ by $(H,\Gamma)$. Let $\Delta_m$ denote the set of $m$-simplices of $L$. Let $\Z\Delta_m$ denote the free abelian group on $\Delta_m$: this is a permutation module for $H$.

We state our results in terms of flag complexes. The stylistic difference from the works \cite{BCK}, \cite{Cornulier} that results from this choice is nothing more than translation using Lemma \ref{cliquestoflags}.

\subsection{Theorem A}\label{s}
\emph{
The following conditions are sufficient for $G$ to be of type $\fn$.
\begin{enumerate}
\item
$H$ is of type $\fn$;
\item
$A$ is of type $\fn$;
\item
$\Z\Delta_p$ is of type $\FP_{n-1-p}$ over $\Z H$ for $0\le p\le n-1$.
\end{enumerate}
}

It may well be possible to establish the natural version of this result for the property $\FP_n$: we conjecture that the result continues to hold when the three instances of $\fn$ are replaced by $\fpn$ in the above statement. Indeed, this is the story presented in \cite{BCK} in the case of wreath products. One can speculate that  in the generality of graph-wreath products, such a generalisation continues to hold, but since it would involve intricate booking keeping with signs in definitions of boundary maps in chain complexes we prefer to leave that case for the present.

In the case of type $\ft$ the sufficient conditions are also necessary.

\subsection{Theorem}\label{n+s2}
\emph{If $A$ is non-trivial then
$G$ is finitely presented if and only if $A$ and $H$ are finitely presented, $\Gamma$ has finitely many orbits of vertices and edges, and each vertex of $\Gamma$ has finitely generated stabilizer.
}

This is easy to prove by the same methods as Cornulier uses for wreath products and we leave the details to the reader.
For completeness we state the very simple result regarding finite generation.

\subsection{Lemma}\label{n+s1}
\emph{If $A$ is non-trivial then
$G$ is finitely generated if and only if both $A$ and $H$ are finitely generated and $\Gamma$ has finitely many orbits of vertices.}
\begin{proof}
Cornulier's proof for wreath products can be employed with essentially no modification, (see \cite{Cornulier}, Proposition 2.1).
\end{proof}

As a step towards finding necessary conditions for the graph-wreath product to have homotopical finiteness we offer the following three results.

\subsection{Theorem B}\label{necessity3+}
\emph{
Suppose that $G$ is of type $\FFF_n$ with $n\ge3$. Assume also that $A$ is non-trivial.
Assume that the stabilizers of $p$-cells of $L$ are of type $\FP_{n-1-p}$ over $\Z H$ for $0\le p\le n-2$.
Then $A$ and $H$ are of type $\fn$ and $H$ acts cocompactly on the $(n-1)$-skeleton of $L$.
}

To obtain necessary and sufficient conditions we make the same assumption that was successfully used in (\cite{BCK}, Theorems A and B).

\subsection{Theorem C}\emph{If $A$ has infinite abelianization
then the following conditions are necessary and sufficient for $G$ to be of type $\FFF_n$.
\begin{enumerate}
\item
$H$ is of type $\FFF_n$
\item
$A$ is of type $\FFF_n$
\item
$\Z\Delta_p$ is of type $\FP_{n-1-p}$ as a $\Z H$-module for $0\le p\le n-1$.
\end{enumerate}
}

Our fourth result does not require the infinite abelianisation condition but instead takes advantage of the special properties of polycylic groups and their groups rings.

\subsection{Theorem D}
\emph{
If $H$ is polycyclic-by-finite and $A$ is non-trivial then $G$ is of type $\FFF_n$ if and only if the following conditions hold:
\begin{enumerate}
\item
$A$ is of type $\FFF_n$.
\item
$H$ acts cocompactly on the $(n-1)$-skeleton of $L$.
\end{enumerate}
}

While it might at first be tempting to conjecture that the assumption that $A_{\mathrm{ab}}$ is infinite in Theorem C can be dropped we advise caution. On the contrary, it seems to us plausible that there are groups $A$ with finite abelianization so that if the clique modules $\Z\Delta_p$ satisfy homological finiteness conditions over finite fields $\Z/p\Z$ for certain $p$ instead of over $\Z$ then the type $\FP_n$ of the graph-wreath product can still be established. Such examples are hard to look for at present because of a lack of examples distinguishing type $\FP_n$ over $\Z$ from type $\FP_n$ over $\Z/p\Z$.

Consider now the special case when $\Gamma$ is the discrete graph on its vertices, that is, $\Gamma$ has no edges; and $A$ is infinite cyclic. Then the base $B$ of the graph-wreath product is the free group on the set $V\Gamma$ of vertices. Thus we have the

\subsection{Corollary}
\emph{
Let $H$ be a group of type $\FFF$ and let $K$ be a subgroup. Let $B$ be the free group on the set of right cosets of $K$ in $H$. Then the semidirect product $B\rtimes H$ is of type $\FFF$ if and only if $K$
is of type $\FP$ over $\Z$.
}

What makes this interesting is that there are examples known in which $H$ is a right-angled Artin group and $K$ is a normal subgroup with $H/K\iso\Z$ such that $K$ is of type $\FP$ but not of type $\FFF$. We refer the reader to \cite{BB} for the details.

\section{Polyhedral Products and the Proofs of Theorems A and B}

Let $L$ be a simplicial complex and let $\mathbf X=(X_v,Y_v,*_v)_{v\in V}$ be a family of pointed pairs of spaces indexed by the vertex set $V$ of $L$. The \emph{finitary product} $\bigoplus_{v\in V}X_v$ is the set of tuples in the full cartesian product $\prod_{v\in V}X_v$ whose supports are finite subsets of $V$. If $L$ is a finite complex then the finitary product coincides with the cartesian product and has the product topology. If $L$ is infinite then the finitary product may be viewed as the colimit of finite products and naturally it is endowed with the colimit topology.
The \emph{grand support} of a tuple $\mathbf x=(x_v)_{v\in V}$ in $\bigoplus_{v\in V}X_v$ is defined by
$$\supp\mathbf x=\{v\in V|\ x_v\notin Y_v\}.$$
The polyhedral product $\mathbf X^L$ is defined to be the subspace of the finitary product comprising tuples whose grand supports belong to $L$. Here we say that a subset $U$ of $V$ belongs to $L$ if it is finite and spans a simplex of $L$. The empty simplex is permitted and so the polyhedral product naturally has the base point
$(*_v)_{v\in V}$. When $L$ is a flag complex with $1$-skeleton $\Gamma$, then a subset $U$ of vertices belongs to $L$ if and only if it is a clique in $\Gamma$.

\subsection{Proof of Theorem A}
Adopt the notation of \S\ref{mainstatement}.
Choose an Eilenberg--Mac Lane space $X$ for $A$ with basepoint $*$ that is witness to $A$ being of type $\FFF_n$. We may choose $X$ to be a CW-complex with a single $0$-cell and with finitely many $p$-cells for $p\le n$. Let $Z$ denote the polyhedral product of the constant family $(X,\{*\},*)_{v\in V}$. Using (\cite{Davis2012}, Theorem 2.22) which is formulated for infinite complexes $L$ as (\cite{DK}, Theorem 1), observe that $Z$ is an Eilenberg--Mac Lane space for the group $B=A^L$ on which the group $H$ acts by self-homeomorphisms. The universal cover $E$ of $Z$ admits an action of the semidirect product $B\rtimes H$, that is an action of $G$. Moreover, the setwise stabilizers of cells in $E$ are isomorphic to subgroups of $H$ that are setwise stabilizers of cliques in the graph $\Gamma$. A cell of dimension $p$ in $E$ may arise as a product of cells from copies of $X$ whose dimensions sum to $p$.

If we assume that there are only finitely many orbits of $(p+1)$-cliques in $\Gamma$ for $p\le n$ then it follows that $E$ has finitely many orbits of $p$-cells for $p\le n$.
In other words, $E$ has finite $n$-skeleton mod $G$.
 If we also assume that the stabilizer in $H$ of each $(p+1)$-clique in $\Gamma$  is of type $\FP_{n-p}$ then it follows that $E$ is $n$-good for the group $G$ in the sense of (\cite{Brown1987}, \S1). It now follows from Brown's elementary criterion (\cite{Brown1987}, Proposition 1.1) that $G$ has type $\FP_n$. When $n\ge2$ we also have that $G$ is finitely presented by Theorem \ref{n+s2}, therefore in all cases, $G$ is of type $\fn$. \qed

\subsection{Proof of Theorem B}

Again, we adopt the notation of \S\ref{mainstatement}.
We now assume that $G$ has type $\fn$. The aim is to investigate what can be said about $H$, $A$, $\Gamma$, and $L$. We may assume that $n\ge3$ since the cases $n\le2$ are understood (see Theorem \ref{n+s2} and Lemma \ref{n+s1}).
It is immediate that $H$ must have type $\fn$ because $H$ is a retract of $G$.

 We work by induction on $n$ and so we may assume that
$A$ is of type $\FFF_{n-1}$, and that $H$ acts cocompactly on the $(n-2)$-skeleton, $L^{n-2}$. Fix an Eilenberg--Mac Lane space $X$ for $A$ with finite $(n-1)$-skeleton.

Using Lemma \ref{cliquestoflags}, the cocompactness of the $H$-action on the $X^{n-2}$-skeleton, and the assumption in our statement of Theorem B we have the following condition.
\begin{itemize}
\item[($\dagger$)]
$\Z\Delta_p$ is of type $\FP_{n-1-p}$ as a $\Z H$-module for $0\le p\le n-2$.
\end{itemize}

To prove that $A$ is of type $\fn$ and that there are finitely many orbits of $(n-1)$-cells in the flag complex $L$ we shall use constructions involving subcomplexes of $L$ and $X$.

\begin{proof}
[Proof that $A$ is of type $\fn$]
Let $(X_\alpha)$ be the family of finite subcomplexes of the Eilenberg--Mac Lane space $X$ which contain the $(n-1)$-skeleton. These are ordered by inclusion and have $X$ as their filtered colimit. Note that since $n\ge3$ these subspaces all have the same fundamental group $A$.

Let $W_\alpha$ be the universal cover of the polyhedral product $X_\alpha^L$. The family of $W_\alpha$ also form a filtered colimit system with colimit $W:=\widetilde{X^L}$. The $W_\alpha$ all have the same $(n-1)$-skeleton as each other and as $W$. Recall that $W$ is an Eilenberg--Mac Lane space for $A^L$. The $W_\alpha$ are therefore $(n-2)$-connected and the structure of their cellular chain complex up to dimension $n-1$ gives the exact sequence
$$0\to Z_{n-1}(W)\to C_{n-1}(W)\to\dots\to C_1(W)\to C_0(W)\to\Z\to0$$
in which $Z_{n-1}(W)$ denotes the group of $(n-1)$-cycles of $W$.
As in the proof of Theorem A the conditions we have through induction and through $(\dagger)$ imply that $C_j(W)$ is of type $\FP_{n-j-1}$ as a $\Z G$-module for $0\le j\le n-1$. The trivial module at the right hand end of this sequence is of type $\fpn$ over $\Z G$ because $G$ is of type $\fpn$. We deduce that $Z_{n-1}(W)$ is finitely generated.
For each $\alpha$ we have short exact sequences
$$B_{n-1}(W_\alpha)\mono Z_{n-1}(W)\epi H_{n-1}(W_\alpha).$$
Since $\colim_{\alpha}H_{n-1}(W_\alpha)=H_{n-1}(W)=0$ and $Z_{n-1}(W)$ is finitely generated we conclude that there is a choice $\alpha_0$ for which $H_{n-1}(W_{\alpha_0})=0$. This tells us that $X_{\alpha_0}$ has trivial $(n-1)$st homology group and since $n\ge3$, the Hurewicz isomorphism holds and $X_{\alpha_0}$ is $(n-1)$-connected. Then an Eilenberg--Mac Lane space for $A$ can be chosen by adding cells of dimension $n+1$ and greater to $X_{\alpha_0}$. Therefore $A$ is of type $\fn$.
\end{proof}

\begin{proof}[Proof that $H$ acts cocompactly on the $(n-1)$-skeleton of $L$.]
For a contradiction, suppose that there are infinitely many orbits of $(n-1)$-simplices in $L$. These orbits are countable in number so we may write the $(n-1)$-skeleton of $L$ as a strictly ascending union of a chain of $H$-finite $H$-complexes beginning with the $(n-2)$-skeleton:
$$L^{n-2}=L_0\subset L_1\subset L_2\subset\dots\subset L^{n-1},$$
$$L^{n-1}=\bigcup_{j\in\N}L_j.$$
Now the spaces $X^{L_j}$, for $j\ge0$ and $X^{L^{n-1}}$ all have the same fundamental group $A^L$ and $X^{L^{n-1}}$ is the union of the
$X^{L_j}$. On passing to universal covers we obtain a chain of $(n-2)$-connected spaces, the $\widetilde{X^{L_j}}$, with union
$\widetilde{X^{L^{n-1}}}$. We may assume that each $L_j$ has been chosen with exactly one orbit of missing $(n-1)$-cells that are included in the next, $L_{j+1}$. We now show that the induced maps in homology:
$$H_{n-1}(\widetilde{X^{L_j}})\to H_{n-1}(\widetilde{X^{L_{j+1}}})$$ have non-trivial kernel for all $j$. Let $\sigma$ be a representative of the orbit of missing $(n-1)$-cells in $L_j$ that are present in $L_{j+1}$. Then $\partial\sigma$ is a subcomplex of $L_j$ and $L_j\cup\sigma$ is a subcomplex of $L_{j+1}$. There is an inclusion and retraction of pairs
$$(X^\sigma,X^{\partial\sigma})\hookrightarrow(X^{L_j\cup\sigma},X^{L_j})\to
(X^\sigma,X^{\partial\sigma}).$$
Thus it suffices to show that the kernel of the map
$$
H_{n-1}(\widetilde{X^{\partial\sigma}})\to H_{n-1}(\widetilde{X^{\sigma}})
$$
is non-zero because we have a commutative square
$$\xymatrix{
H_{n-1}(\widetilde{X^{\partial\sigma}})\ar[r]\ar[d]
&H_{n-1}(\widetilde{X^{\sigma}})\ar[d]\\
H_{n-1}(\widetilde{X^{L_j}})\ar[r]
&H_{n-1}(\widetilde{X^{L_{j+1}}})
}$$
in which the vertical maps are injective.
Note that $X^{\partial\sigma}$ and $X^\sigma$ have the same $(n-1)$-skeleton: $\sigma$ is the flag complex on the complete graph on its vertices, so 
$\widetilde{X^{\sigma}}=\widetilde X^\sigma$ and since
$n\ge3$, there is an inclusion of $\widetilde X^{\partial\sigma}$ into $\widetilde X^{\sigma}$ because $X^{\partial\sigma}$ and $X^\sigma$ both have the same fundamental group. Moreover, $\widetilde X^{\partial\sigma}$ and $\widetilde{X^{\partial\sigma}}$ are the same since they are both simply connected. So it suffices to prove that the map
$$
H_{n-1}(\widetilde X^{\partial\sigma})\to H_{n-1}(\widetilde X^{\sigma})
$$
has non-trivial kernel.
Just as in the proof that $A$ is of type $\FFF_n$, the $(n-1)$-cycles of the complexes $\widetilde{X^{L_j}}$ are independent of $j$ and form a finitely generated $\Z H$-module. We reach a contradiction by essentially the same reasoning as in the proof that $A$ is of type $\fn$.
\end{proof}

\section{A little homological algebra and the Proof of Theorem C}

Let $J$ be a commutative ring. Usually, $J$ is $\Z$ or a prime field.
We shall need the next result only in case $J=\Z$. However the result is inspired by classical results of Nakaoka and a useful account can be found in (\cite{Evens}, Chapter 5) where a similar argument is employed to describe the mod $p$ cohomology rings of Sylow $p$-subgroups of certain finite symmetric groups. Here, we shall use the result as a key step in the proof of Theorem C. The formulation of this Proposition has been influenced by discussions between the first author and David Benson. Originally this reasoning was inspired by (\cite{BCK}, Lemma 4.1 and Proposition 4.2). The statement involves a group $G$ with a normal subgroup $B$ and we are interested in resolutions $P_*\epi J$ which are on the one hand exact sequences of $JG$-modules and $JG$-maps
$$\dots\to P_n\to P_{n-1}\to\dots\to P_1\to P_0\to J\to 0$$
while on the other hand are $JB$-projective, meaning that each $P_j$ is projective as a $JB$-module. We refer to such exact sequences and $JB$-projective $JG$-resolutions.

\subsection{Proposition}\label{nakaoka}
\emph{
Let $G$ be a group with a normal subgroup $B$. Suppose that there is a $JB$-projective $JG$-resolution $P_*\epi J$ of the trivial module such that
\begin{enumerate}
\item[$(*)$]
the induced maps $P_j\otimes_{JB}J\to P_{j-1}\otimes_{JB}J$ are zero for all $j\ge1$.
\end{enumerate}
Then, for all $p\ge0$, we have natural isomorphisms as follows.
\begin{enumerate}
\item
If $M$ is any left $JG$-module on which $B$ acts trivially
$$H_p(G,M)\iso\bigoplus_{i+j=p}\tor_i^{JH}(H_j(B,J),M).$$
\item
If $N$ is any right $JG$-module on which $B$ acts trivially
$$H^p(G,N)\iso\bigoplus_{i+j=p}\ext^i_{JH}(H_j(B,J),N).$$
for any right $JG$-module $N$.
\end{enumerate}
}
\begin{proof} The hypothesis $(*)$ implies that $H_j(B,J)=P_j\otimes_{JB}J$, and this is a free $J$-module for all $j$.
Thus, using the universal coefficient theorem, we have $$H_j(B,M)=H_j(B,J)\otimes_JM,$$ when $M$ is a left $JG$-module on which $B$ acts trivially.
Write $H$ for the quotient $G/B$.
Whether or not $B$ acts trivially on $M$, there is a Lyndon--Hochschild--Serre spectral sequence
$$E^2_{ij}=H_i(H,H_j(B,M))\implies H_{i+j}(G,M).$$

Let $M$ be any left $JH$-module. Note that $M$ can also be viewed as a $JG$-module on which $H$ acts trivially. The spectral sequence takes this form:
$$E^2_{ij}=\tor_i^{JH}(H_j(B,J),M)\implies H_{i+j}(G,M).$$
The hypothesis $(*)$ can be used to show that the differentials $d^\ell_{ij}$ all vanish for $\ell\ge2$ and this collapsing of the spectral sequence leads to the desired conclusion (i). To see why this happens is not hard so we include the details by manufacturing the spectral sequence from scratch: choose any resolution $Q_*\epi M$ by $JH$-projective left $JH$-modules and consider the double complex
$$E^0_{ij}=P_j\otimes_JQ_i.$$
Here, $G$ acts diagonally on $P_j\otimes_JQ_i$ by $g(\xi\otimes\eta)=\xi g^{-1}\otimes g\eta$, the total complex $E^{tot}_*$ resolves $M$ by $JG$-projective left $JG$-modules, and the homology of the total complex yields $H_*(G,M)$. That is, the homology groups $H_*(G,M)$ are the homology groups of the chain complex $J\otimes_{JG}E^{tot}_*$ where $E^{tot}_*$ denotes the total complex
$$E^{tot}_p=\bigoplus_{i+j=p}P_j\otimes_JQ_i.$$
The spectral sequence viewpoint comes from writing the functor
$J\otimes_{JG}$ as the composite $J\otimes_{JB}$ followed by
$J\otimes_{JH}$. Applying the hypothesis $(*)$ again and also using that $B$ acts trivially on the $Q_i$ we have
$$J\otimes_{JG}E^{tot}_p
=J\otimes_{JH}(J\otimes_{JB}E^{tot}_p)
=J\otimes_{JH}(\bigoplus_{i+j=p}H_j(B,J)\otimes_JQ_i).$$
But now the vertical maps
($H_j(B,J)\otimes_JQ_i\to H_{j-1}(B,J)\otimes_JQ_i$) in the intermediate double complex
$$J\otimes_{JB}E^{tot}_*$$
 are all zero so we at once see that its $p$th homology is
$$\bigoplus_{i+j=p}\tor_i^{JH}(H_j(B,J),M).$$
This completes the argument for (i). The argument for (ii) is similar.
\end{proof}

\subsection{Corollary}
\emph{
Let $B\normal G$ be groups satisfying hypothesis $(*)$ of Proposition \ref{nakaoka}. If $G$ is of type $\FP_n$ over $J$ then $H_p(B,J)$ is of type $FP_{n-p}$ as a $JH$-module for $0\le p\le n$. In particular, $H:=G/B$ is of type $FP_n$ over $J$.
}

\begin{proof}
According to (\cite{Bieriqmw}, Theorem 1.3 (i)$\iff$(iiib)) a module $M$ over an associative ring $R$ is of type $\FP_m$ if and only if the functor $\ext_R^i(M,\quad)$ commutes with filtered colimits for each $i<m$. Thus, if $G$ is of type $FP_n$ over $J$ then $H^{p}(G,\quad)$ commutes with filtered colimits for all $p<m$. Proposition \ref{nakaoka}(ii) shows that there is an isomorphism of functors of $JH$-modules:
$$H^{p}(G,\quad)\iso\bigoplus_{i+j=p}\ext^i_{JH}(H_j(B,J),\quad),$$
and we infer that all the functors
$$\ext^i_{JH}(H_j(B,J),\quad)$$ commute with filtered colimits when $i<j<m$. Thus, for each $j$, $H_j(B,J)$ is of type $\FP_{m-j}$.

Alternatively one can also reason using Proposition \ref{nakaoka}(i) together with the criterion that the module $M$ over $R$ is of type $\FP_m$ if and only if $\tor^R_i(M,\prod R)$ vanishes for $0<i<m$ and there is a natural isomorphism 
$\tor^R_0(M,\prod R)=\prod\tor_0^R(M,R)$, (see \cite{Bieriqmw}, Theorem 1.3, (i)$\iff$(iiia). Arguably this second approach is simpler save for the issue that the case $i=0$ must be handled differently from the cases $i>0$.
\end{proof}

\subsection{Proof of Theorem C} Suppose now that, with the same notation as in \S\ref{mainstatement}, the group $G$ is of type $\FFF_n$. We use the hypothesis that $A/[A,A]$ is infinite to deduce properties (i), (ii), (iii) of Theorem C. In the matter at hand, $A$ is always a finitely generated group and hence there is a surjective homomorphism $A\to\Z$ witnessing $\Z$ as a retract of $A$. The graph-wreath product construction behaves functorially and so we deduce that $\Z\wr_\Gamma H$ is a retract of $A\wr_\Gamma H$. The property of being of type $\fn$ is inherited by
$\Z\wr_\Gamma H$. Therefore we may assume that $A=\Z$. From the fact that $G$ is of type $\fpn$ we infer, using Proposition \ref{nakaoka}
and its Corollary from this section, that $H_p(B,\Z)$ is a module of type $\FP_{n-p}$ for $0\le p\le n$. Now $B$ is a right-angled Artin group and the homology group $H_p(B,\Z)$ is the free abelian group on the $(p+1)$-cliques in the graph $\Gamma$. This puts in place the hypothesis  $(*)$ that is required in the Proposition and Corollary. We refer the reader to Charney's survey \cite{Charney} for this fact and background information about right angled Artin groups. The result follows from Theorem B.

\section{The proof of Theorem D}

If $n=1$ then Theorem D follows immediately from Lemma \ref{n+s1}.
If $n=2$ then it is the special case of Theorem \ref{n+s2} when $H$ is polycyclic-by-finite.  From now on we assume that $n\ge3$ and proceed by induction.

If $G$ is of type $\FFF_n$ then by induction we may assume that $H$ acts cocompactly on the $(n-2)$-skeleton of $L$ and also that $A$ is of type $\FFF_{n-1}$. In particular, $A$ is finitely presented and therefore, by Theorem \ref{n+s2}, $\Z\Delta_0$ is finitely presented and $\Z\Delta_1$ is finitely generated over $\Z H$. But $\Z H$ is a Noetherian ring because $H$ is polycyclic-by-finite (see (\cite{Hall} Theorem 1 and the immediately following paragraph) or (\cite{LR}, 4.2.3)) , and therefore 
$\Z\Delta_1$ and $\Z\Delta_2$ are modules of type $\fpinfty$ over $\Z H$. By induction we may assume that $\Z\Delta_{i}$ is finitely generated (and therefore $\fpinfty$) for $i<n-1$.
This ensures that the hupotheses of Theorem B hold and it follows that $A$ is of type $\FFF_n$ and that $H$ acts cocompactly on the $(n-1)$-skeleton of $L$. 

The converse, that conditions (i) and (ii) of Theorem D imply $G$ to be of type $\FFF_n$ is a direct consequence of Theorem A together with the information that $H$ is polycyclic-by-finite. Again, we may appeal to the fact that $\Z H$ is Noetherian so that $H$ is of type $\FFF_\infty$ and all finitely generated $\Z H$-modules are of type $\fpinfty$. 

This completes the proof of Theorem D. Notice that the hypothesis that $A$ has infinite abelianisation is not required. 

\section{Closing Remarks}

\subsection{Houghton Groups}

For $n\ge1$, let $H_n$ denote the group of those permutations of the set $R_n:=\{1,\dots,n\}\times\N$ which act as translation far along any individual ray $\{i\}\times\N$, $1\le i\le n$. Brown \cite{Brown1987} proved that $H_n$ is of type $\FFF_{n-1}$ but not of type $\fn$ for each $n\ge1$. The pointwise stabilizer in $H_n$ of any finite subset of $R_n$ is easily seen to be isomorphic to $H_n$ and the setwise stabilizer is a finite extension of the pointwise stabilizer: therefore all these stabilizers have type $\FFF_{n-1}$.
 Moreover, for all $p$, $H_n$ acts transitively on sets of cardinality $p$. Thus we have the following consequence.

\subsection{Corollary}
\emph{
A wreath product $A\wreath_{R_n} H_n$ has type $\FFF_{n-1}$ if and only if $A$ is of type $\FFF_{n-1}$.
}

By taking iterated wreath products of Houghton's groups we may state:

\subsection{Corollary}
\emph{
For any $m$ and $k$ there exists elementary amenable groups of type $\FP_m$ whose F{\o}lner set growth functions grow at least as fast as $n^{n^k}$.
}

\begin{proof}
Houghton's groups are locally finite by abelian and hence are amenable. They belong to the class of elementary amenable groups: the smallest class containing all abelian and finite groups which is also closed under extensions and directed unions.
Using the observation that for $n\ge2$, Houghton's group $H_n$ contains an infinite cyclic group $C$ which has a regular orbit (i.e. a $C$-orbit on which the $C$ acts freely) on the ray set $R_n$ we may infer that the iterated wreath product
$$(\dots(H_{m+1}\wreath_{R_{m+1}}H_{m+1})\wreath_{R_{m+1}}H_{m+1}\dots)\wreath_{R_{m+1}}H_{m+1}$$ contains subgroups isomorphic to the iterated standard wreath product
$$(\dots(\Z\wreath\Z)\wreath\Z\dots)\wreath\Z.$$
Now our result follows from calculations of Erschler (\cite{Erschler}, Theorem 1 and Example 3).
\end{proof}

Another interesting source of examples arise from Thompson's groups. The homological finiteness properties of these groups is established by Brown \cite{Brown1987} and we refer the reader to that paper for background material. For example, Thompson's group $F$ acts on the set $D:=\Z[\frac12]\cap[0,1]$ of dyadic rationals in the unit interval with finite set stabilizers that are finite products of copies of Thompson's group and therefore are all groups of type $\fpinfty$. Therefore we may state:

\subsection{Corollary}
\emph{
The wreath product $A\wreath_DF$ is of type $\FP_n$ whenever $A$ is of type $\FP_n$.
}

Taking $A$ to be first a non-trivial finite group, and secondly an infinite cyclic group, we have corollaries as promised in our introduction.

\subsection{Corollary}
\emph{
There exist groups of type $\fpinfty$ with infinite locally finite normal subgroups. There exist groups of type $\fpinfty$ with infinite rank free abelian normal subgroups.
}

\section{Acknowledgements}

The first author wishes to acknowledge that his discussions at various times with David Benson, Michael Davis, Ross Geoghegan, Dessislava Kochloukova, and Ian Leary have led to considerable improvements to the results and proofs here recorded.


\end{document}